\documentclass{article}
\usepackage{latexsym}
\usepackage{amssymb}

\advance \textwidth by 2.1cm
\advance \textheight by 3.6cm
\advance \topmargin by -1.3cm
\advance \oddsidemargin by -1.2cm

\newfont{\klein}{cmr5}
\newfont{\symklein}{cmsy5}
\newfont{\mathsymfett}{cmbsy10}
\newfont{\sansserifklein}{cmss8}

\newtheorem{lemma}{Lemma}
\newcommand{\be}{\begin{equation}}
\newcommand{\ee}{\end{equation}}

\newcommand{\EinrL}[1]{\hangindent=1.85cm\hangafter=1 #1}

\newcommand{\sst}{\scriptscriptstyle}
\newcommand{\s}[1]{\mbox{$\scriptscriptstyle{#1}$}}
\newcommand{\ssf}[1]{\mbox{$\mbox{\sansserifklein #1}$}}
\newcommand{\bom}[1]{\mbox{\boldmath $#1$\unboldmath}}

\newcommand{\mde}{\ensuremath{\, \mid \, }}

\newcommand{\Rl}{\ensuremath{\mathbb{R}}}
\newcommand{\Cx}{\ensuremath{\mathbb{C}}} 
\newcommand{\Nl}{\ensuremath{\mathbb{N}}}
\newcommand{\Cinfty}{\ensuremath{C \rule{0ex}{1.8ex}^{\vspace{-0.15ex} \scriptscriptstyle \infty}}}
\newcommand{\CP}{\mbox{$C \rule{0ex}{1.8ex}^{\vspace{-0.15ex} \scriptscriptstyle \infty}(P)$}}  
\newcommand{\CPl}{\ensuremath{C \rule{0ex}{1.8ex}^{\vspace{-0.15ex} \scriptscriptstyle \infty}(P)
               [ \! [ \lambda ] \! ] }} 
\newcommand{\CUl}{\ensuremath{C \rule{0ex}{1.8ex}^{\vspace{-0.15ex} \scriptscriptstyle \infty}(U)
               [ \! [ \lambda ] \! ] }}          
  
\newcommand{\CCdo}{\mbox{$C \rule{0ex}{1.8ex}^{\vspace{-0.15ex} \scriptscriptstyle \infty}
               (\hat{\cal C},\hat{\omega})$}}      
\newcommand{\CCl}{\ensuremath{C \rule{0ex}{1.8ex}^{\vspace{-0.15ex} \scriptscriptstyle \infty}
               ({\cal C})[ \! [ \lambda ] \! ] }} 
\newcommand{\CCdl}{\ensuremath{C \rule{0ex}{1.8ex}^{\vspace{-0.15ex} \scriptscriptstyle \infty}
               (\hat{\cal C})[ \! [ \lambda ] \! ] }} 
\newcommand{\CPn}{\ensuremath{\mathbb{C}P^n}}               
\newcommand{\Cxwz}{\mbox{$\mathbb{C}^{\,n}\rule{0ex}{1.6ex}^{\!\sst +}{}^1
               \!\setminus\!\{0\}$}}

\newcommand{\imbed}{\ensuremath{\hookrightarrow}}

\newcommand{\lra}{\mbox{$\longrightarrow$}}

\newcommand{\BBox}{\vrule width 1.6mm height 1.6mm depth 0mm}

\newcommand{\halb}{\mbox{$\frac{1}{2}\,$}}
\newcommand{\zzq}{z\bar{z}}
\newcommand{\pr}{\ensuremath{[\![ \lambda ]\!]}}

\newcommand{\la}{\ensuremath{\lambda}}
\newcommand{\Cc}{{\ensuremath{\cal C}}}
\newcommand{\Cd}{{\ensuremath{\hat{\cal C}}}}
\newcommand{\Eq}{{\ensuremath{\bar{E}{}}}}

\newcommand{\id}{\ensuremath{\mbox{\sf Id}}}

\newcommand{\prol}{\ensuremath{\mbox{\sf prol} }}
\newcommand{\supp}[1]{\ensuremath{\mbox{\sf supp}(#1)}}
\newcommand{\p}{\ensuremath{\mbox{\sf p}}}
\newcommand{\pj}{\ensuremath{\pi_{\!\scriptscriptstyle J}} }               

\newcommand{\Mn}{\ensuremath{\rule{0ex}{2.3ex}^{\scriptscriptstyle 
               \longleftrightarrow} \hspace*{-3.1ex}M_n}}
\newcommand{\M}[1]{\ensuremath{\rule{0ex}{2.3ex}^{\scriptscriptstyle 
               \longleftrightarrow} \hspace*{-3.1ex}M_{#1}}}
\newcommand{\resC}{\ensuremath{\mbox{\sf res}}_{\cal C}}

\title {Star Product Reduction\\
        for Coisotropic Submanifolds of Codimension 1}

\author {{\bf Peter
          Gl\"o\ss{}ner\thanks{pegl@phyq1.physik.uni-freiburg.de}
         } \\[3mm]
         Fakult\"at f\"ur Physik\\Universit\"at Freiburg \\
         Hermann-Herder-Str. 3 \\
         79104 Freiburg i.~Br., F.~R.~G \\[3mm]
        }

\date{FR-THEP-98/10 \\[1mm]
      20 April 1998 \\[5mm]}

\begin{document}

\maketitle

\begin {abstract}
We propose a reduction procedure that leads to a reduced star product on the reduced 
phase space of a ``First Class''--constrained system, where no symmetries, group actions 
or the like are present. For the case that the coisotropic constraint submanifold has 
codimension 1, we establish a constructive method to compute the reduced star product 
explicitly. Concluding examples show that this method depends crucially on the constraint 
function singled out to describe the constraint submanifold and not only on this 
submanifold itself, and that two different constraint 
functions for the same constraint submanifold will generally result in not only different 
but inequivalent reduced star products.
\end {abstract}

\subsection{Introduction} \label{Intro}

In 1978, Bayen, Flato, Fr{\o}nsdal, Lichnerowicz and Sternheimer
established in \cite{BFFLS78} the concept that we call {\bf deformation 
quantization}. The idea is to replace the pointwise product on the 
algebra of phase space functions $\CP$ of a classical dynamical system 
by a non--commutative so--called {\bf star product} $\ast$ on $\CPl$.
For two $f$, $g \in \CPl$, $f \ast g$ is again a formal power series 
in $\la$, is the pointwise product in $0$th order, and $f \ast g - g
\ast f$ is $i$ times the Poisson bracket $\{ f,g \}$. In this way, the
formal parameter $\la$ plays the r\^{o}le of Planck's constant $\hbar$
and the algebra $\CP$ can be addressed as the algebra of quantum 
observables, deformed from the classical observables $\CP$ in the sense
of Gerstenhaber \cite{GS88}. The naturally arising question whether star
products exist for arbitrary symplectic manifolds was answered in the 
affirmative 1983 by DeWilde and Lecomte \cite{DL83}, Omori, Maeda and
Yoshioka \cite{OMY91}, while Fedosov \cite{Fed94} gave a more geometric proof 
in 1994. In the more general setting of a Poisson manifold, the existence 
of star products was shown only recently by Kontsevich \cite{Kon97}.

In this paper, we propose a kind of ``quantum'' or ``deformed'' reduction 
mechanism. It shall serve to answer the following question: given an arbitrary {\bf star 
pro\-duct} $\ast$ on $(P,\omega)$, and a coisotropic submanifold $\Cc \hookrightarrow P$ 
with Reduced Phase Space $(\hat{\Cc}, \hat{\omega})$, how can we construct a star 
product $\star$ on the symplectic manifold $(\hat{\Cc}, \hat{\omega})$ which can be 
addressed as the result of some reduction procedure starting with $(P, \ast)$?

This problem has been dealt with in several publications: after the authors
of \cite{BFFLS78} had discussed some basic examples, Fedosov \cite{Fed94b} 
investigated reductions induced by certain $U(1)$--actions on the phase space.
The case of $\CPn$ and its noncompact dual were considered in \cite{BBEW96a}
and \cite{BBEW96b} by Bordemann et al., even resulting in explicit formulae. 
Schirmer \cite{Sch97} gives a generalization to Grassmann manifolds. In 
\cite{Fed97}, Fedosov introduces a reduction procedure for Hamiltonian group
actions of arbitrary compact Lie groups. 

But in these treatments it has been 
essential that the constraint submanifold $\Cc$ was a level surface of a momentum mapping,
that is, the phase space reduction in these cases was the outcome of some symmetry group 
acting on the phase space $P$; all reduction processes cited above consider 
a star product on the phase space $P$ such that the invariant functions on 
$P$ form a subalgebra.
We will try a first step towards an algebraic reduction 
process for general coisotropic submanifolds, as they may occur for example as the result 
of a non--surjective fiber derivative (see \cite{AM85}, esp. 3.5).
Symmetry groups, momentum mappings etc.\ will 
consequently play no part in our considerations, and we are free to choose
any $f \in \CP$ vanishing on $\Cc$ as constraint function. 

Before we sketch our methods and results, we briefly review the notion of
classical phase space reduction , in order to introduce the algebraic structures 
we will deform in the ensuing section. As general reference for this subject
may serve \cite{AM85}, Chapter 5, esp. exercises and references therein.
\newline If a dynamical system on a phase space $P$ (equipped with symplectic form $\omega$) is 
forced to stay on a coisotropic  submanifold $\Cc \imbed P$ (``by first class 
constraints''; see \cite{D64} and \cite{GNH78}), we can perform the well-known procedure 
commonly called {\bf phase space reduction} \cite{AM85} to obtain a phase space which represents 
in some sense the true degrees of freedom of the physical system. The reduction process in its 
differential geometric picture consists essentially in pulling back $\omega$ on $\Cc$ and 
then dividing $\Cc$ by the foliation which is generated by the integrable distribution 
associated to the kernel of the pull--back of $\omega$. We come out with the {\bf Reduced 
Phase Space} $(\hat{\Cc}, \hat{\omega})$. On the algebraic side, this picture is reflected 
in the following way. The constaint manifold $\Cc$ is characterized by its vanishing ideal 
$I := 
\{ f \in \CP \mde i^{*}f = 0 \}$ with $i: \Cc \hookrightarrow P$ the embedding. This 
algebra is a Poisson ideal in $B := \{ f \in \CP \mde \{f,f'\} \in I 
\:\:\forall \: f' \in I \}$ ($\{ \cdot , \cdot \}$ denoting as usual the Poisson bracket that 
comes with $\omega$), which in turn is a Poisson subalgebra of $\CP$. We therefore can 
define the quotient $B/I$ as the {\bf Reduced Algebra} of the constrained system. It 
carries a Poisson structure inherited from that of $P$, and it turns out that $B/I$ is 
Poisson--isomorphic to $\CCdo$, the functions on the Reduced Phase Space (this is the case 
essentially because the Hamiltonian vector fields to functions in $I$ span the kernel of 
the pull--back of $\omega$). In this sense, the construction of $B/I$ is the algebraic 
form of phase space reduction. 

In the next section, we propose a deformation of the classical algebras $I$ 
and $B$ into new algebras $I^{\ast}$ and $B^{\ast}$. The quotient algebra
$B^{\ast}/I^{\ast}$ turns out to be an associative star algebra, so that
the aim is to establish a linear isomorphism from $B^{\ast}/I^{\ast}$ to
$B/I$ which then is declared to be a star product isomorphism, thereby 
providing $B/I$ with the desired star product addressed as ``reduced from
$(P,\ast)$''. We emphasize that -- though the algebras $B/I$ 
and $B^{\ast}/I^{\ast}$ can always be formed -- the construction of such 
an isomorphism is (of course) by no means natural and constitutes the 
essential task of our reduction process.
\newline Sections 3, 4 and 5 concentrate on the codimension $1$ case and provide
the necessary structures and proofs for a constructive method to compute
this isomorphism. 
\newline Sections 6, 7 and 8 give examples to show plausibility and
feasibility of the reduction process. A simple $\Rl^{2n}$ example is 
followed by a reduction of the Wick product to $\CPn$, where the reduced
star product is the same as in \cite{BBEW96a}, although constructed in a 
completely different way. But, as mentioned before, we can choose our
constraint function arbitrarily among all functions in $I$ and are in no
way restricted to stick to the $U(1)$ momentum mapping that usually serves to the 
classical reduction $\Cxwz \lra \CPn$. We make use of this our freedom in section 8
and reduce the Wick product again, this time by a different constraint 
function. The obtained reduced star product is different from the one derived
before, and what is more, it is even inequivalent to this, as our concluding
remarks will show. This may in turn be set in contrast to \cite{W98}, 
where an inequivalent reduced star product on $\CPn$ is obtained by deforming
the classical momentum mapping, that is by adding terms of higher order in 
the formal parameter but leaving the classical momentum mapping untouched
in $0$th order.

But before we start our discussions, we fix our notion of a star product $\ast$:
a star product of two formal power series $f, g \in \CPl$ be defined as $f 
\ast g := \sum_{\sst n = 0}^{\sst \infty} \la^n M_n(f,g)$, where $M_0(f,g) = f 
\cdot g$ is the usual point product, $M_1(f,g) - M_1(g,f) = i\{ f,g \}$ is $i$ times the 
Poisson bracket on $P$, $f \ast 1 = 1 \ast f = f$ and $\supp{f \ast g} \subseteq \supp{f}
\,\cap\, \supp{g}$, the latter condition saying that the star product is local. We will 
sometimes need the antisymmetric parts of the $M_n$'s (or rather two times this) and 
denote them by $\Mn(f,g) := M_n(f,g) - M_n(g,f)$.

\subsection{Translating the classical into deformed structures}

To establish a star product phase space reduction, we define algebraic structures 
corresponding to the classical ones $I$ and $B$ reviewed in section \ref{Intro}. (We 
remark that from now on, $I$ and $B$ shall be considered as containing {\bf formal power 
series} from $\CPl$ instead of simple functions from $\CP$.) We begin with $I$, which we 
want to ``deform''  into an $I^*$. What properties should this have? 
\newline Firstly, since in a theory quantized by representation of the observables as 
operators,
the Hilbert space of physical states may often be defined via the operators $\hat{J}_i$ 
corresponding to the constraint functions $J_i$, $i = 1 
\ldots \mbox{\sl codim}\: \Cc$, we make the ansatz that such a star product reduction  
incorporates a preferred choice of codimension--$\Cc$--many first class 
constraints. This way, the star product reduction ``sees'' not only the constraint surface 
$\Cc$ itself, but a certain sandwich neighborhood $U \supset \Cc$, as different from the 
classical phase space reduction. So, we construct the new algebra $I^*$ to contain the first 
class constraints $J_i$, $i = 1 \ldots \mbox{\sl codim}\: \Cc$ in an explicit way. We take 
a special set of constraint functions as given and refrain from discussing the reasons 
that could lend preference to this choice over possible other ones; the reasons may be 
found in symmetries of $(P, \omega)$, or there may be no reasons at all -- the star 
product reduction should work with every set (but dependent on it). 
\newline Secondly, $I^*$ should contain the (noncommutative) star product and therefore 
be a one--sided star ideal. 
\newline Thirdly, the classical vanishing ideal $I$ should be regained by performing the 
limit $\la \longrightarrow 0$, $\la$ being the formal deformation parameter. 
\newline These considerations lead to the definition of the {\bf star--left--ideal}
     \be \textstyle
          I^* := \{ f \in \CPl \mde f = \sum g^i \ast J_i \:\:\mbox{\it  for 
          some }\: g^i \in \CPl\,; \:i = 1 \ldots \mbox{\it codim}\, \Cc \}.                            
     \ee 
     
We proceed along these lines defining $B^*$ as ``deformation'' of $B$. $B^*$ has to be a 
star--subalgebra of $\CPl$, and $I^*$ has to be a two--sided star--ideal in $B^*$.
     \be 
          B^* := \{ f \in \CPl \mde f\ast g - g\ast f \in I^* \:\: \forall \: g \in I^* \} 
     \ee                                                                                        
fulfills these requirements, as can be seen by simple computations. So, $B^*/I^*$ is 
well--defined and carries by representant--wise definition $[f] \ast [g] := [f \ast g]$ a 
well--defined star product, thereby forming an associative star algebra. We will use $B^*$ 
mostly in the equivalent form of 
     \be \label{B^*}
          B^* = \{ f \in \CPl \mde J_i \ast f \in I^* \:\: \forall \: J_i ,                  
          i = 1 \ldots \mbox{\it codim}\: \Cc \} \,\,,                                        
     \ee                                                                                   
which can be obtained from the definition in a rather direct way. 

The general aim now is to establish a linear isomorphism between $B/I$ (enlarged, as we 
agreed, to obtain formal power series), and $B^*/I^*$. Then the isomorphisms
     \begin{displaymath}
          \CCl \:\cong\: B/I \:\cong\: (B^*/I^*, \ast)                                                  
     \end{displaymath}                                                                           
provides us with the desired reduced star product on the classical Reduced Phase Space. 
\newline In the next three sections, we will construct the isomorphism $B/I \cong B^*/I^*$ 
in an explicit way, but only for the codimension $1$ case.

\subsection{Sum decomposition of $C^{{}^{\scriptscriptstyle}\infty}(P)[ \! [ \lambda ] \! ] $}

Our first step consists in defining a prolongation prescription for series of functions on 
$\Cc$, i.\ e.\ a mapping 
     \begin{eqnarray*}
           & & \mbox{\sf p} : \CCl \longrightarrow \CUl 
           \hspace{1.5ex} \mbox{with} \hspace{1.5ex} i^{\ast} \circ \mbox{\sf p} = \id \\                                           
           & & (\, U \subset P \:\mbox{\it is a sandwich neighborhood of the 
           constraint surface}\: \Cc \,) \,\,.                                                                                                                                            
     \end{eqnarray*} 
This can be done arbitrarily; we could, for instance, establish a Riemannian metric $g$ on $P$ 
and use the gradient flow of the constraint function -- remember that we are in the 
codimension $1$ case. In many examples with symmetry, a preferred choice for this 
prescription will present itself. However, the reduction process works for every choice 
(but is dependent on it).  
\newline Let $I : \Cc \imbed P$ be the imbedding of the constraint surface into $P$, we 
then define the {\bf prolongation} of a series $f \in \CPl$ by
     \begin{eqnarray}
          & & \prol : \CPl \longrightarrow \CUl \nonumber \\                                
          & & \prol (f) := \mbox{\sf p} (i^*f) \:,                                                       
\end{eqnarray} and we set 
     \be
          F := \{ f \in \CPl \mde f(p) = (\prol (f))(p) \:\: \forall \: p \in U \} \:,
     \ee calling such series ``pure prolongations''. We agree that from now on, we do not 
distinguish between (series of) functions that are different just outside a sandwich 
neighborhood $U$. ``Uniqueness'' will be understood in this sense in what follows.---  Now, 
because clearly $f - \prol(f) \in I$ for all $f \in \CPl$, ``Hadamard's trick'' or any 
other form of the mean value theorem gives us a unique smooth series $h \in \CUl$ such 
that $f - \prol(f) = h \cdot J$. We set $\pj(f) := h$ as the ``component of $f$ along the 
constraint function $J$''. Remark that while $\prol$ is a projection, $\pj$ is not. We end  
with a uniquely defined decomposition of $\CPl$ as a direct sum
     \begin{equation} \label{pijott}
     \begin{array}[b]{rcccl} 
          \CPl & = &        F & \oplus & I \\
             f & = & \prol(f) &      + & \pj(f) \cdot J
     \end{array}           
     \end{equation}  
     
In a second step, we inductively define the following formal power series of $\Cx$--linear 
operators on $\CPl$, using the bilinear operators $M_r$ of our star product $\ast$ on $P$,
the $\pj$ just introduced and the constraint function $J$:
     \begin{eqnarray} \label{Tn} 
          T &=& \textstyle \sum_{n=0}^{\infty} \la^n T_n \nonumber \\
             & & T_0 := \id \nonumber \\
             & & \textstyle T_n(f) := -\sum_{k=1}^{n} T_{n-k}(M_k(\pj(f),J)) 
                 \:\:\mbox{\it for}\:\: n \ge 1
     \end{eqnarray}
     
\begin{lemma} \label{T} \EinrL The above defined operator $T$ has the following properties:
\newline i) $T : I^* \rightarrow I$ is one--to--one and onto.
\newline ii) $T(\prol(f)) = \prol(f)$ for any $f \in \CPl$, and $T(J) = J$.
\newline iii) $T(f \ast J) = f \cdot J$ for any $f \in \CPl$.
\newline iv) for every $g \in \CP$, $\supp{T_n(g)} \cap \Cc \:\subseteq\: \supp{g} \cap \Cc$
             for all $n$.
\end{lemma}
   
{\bf Proof.} i) Injectivity follows from $T_0 = \id$ and surjectivity is clear once iii) 
is proven. ii) follows from the facts that $\pj \circ \prol = 0$ and $\pj(J) = 1$. iii) It 
is sufficient to consider {\it functions} $f \in \CP$. On one hand it is in $n$th order 
$(T(f \ast J))_n = T_n(fJ) + \sum_{k=1}^{n} T_{n-k}(M_k(f,J))$, $n \geq 1$. On the other 
hand, $T_n(fJ) = -\sum_{k=1}^{n} T_{n-k}(M_k(\pj(fJ),J))$ -- but $\pj(fJ) = f$, so $(T(f 
\ast J))_n = 0$ for $n \geq 1$. $(T(f \ast J))_0 = fJ$ is trivial. iv) We use $\supp{\pj(g)} 
\cap \Cc \:\subseteq\: \supp{g} \cap \Cc$ and the fact that the $M_k$ do not enlarge the 
supports of their arguments in a straightforward inductive reasoning. \BBox

This operator $T$ yields another sum decomposition of $\CPl$. Indeed, if we apply $T^{-1}$ 
to $T(f) = \prol(T(f)) + \pj(T(f)) \cdot J$, keeping in mind properties ii) and iii) from 
the above lemma, we obtain for every $f \in \CPl$ the equation $f = \prol(T(f)) + 
\pj(T(f)) \ast J$, the last summand being an element of $I^*$. So there is a further sum 
decomposition of $\CPl$ as 
     \begin{equation} \label{decompT}
     \begin{array}[b]{rcccl} 
          \CPl & = &          F  & \oplus & I^* \\
             f & = & \prol(T(f)) &      + & \pj(T(f)) \ast J 
     \end{array}           
     \end{equation}

\subsection{Isomorphism between $B \cap F$ and $B^* \cap F$}

With the help of the operator series $T$ just defined, we are able to construct a formal 
power series of $\Cx$--linear operators $S_n$ which establishes a linear isomorphism $S : 
B \cap F \rightarrow B^* \cap F$, needed to map $B/I$ and $B^*/I^*$ on each other 
$\Cx$--linearly and bijectively. 
We point out that we will now make an additional assumption, namely, we suppose
there is a transversal section $\sigma$ of the foliation on $\Cc$ associated to 
the Hamiltonian vector field $X_J$ of the constraint function $J$, and with
${\mathfrak p} : \Cc \rightarrow \sigma$ we denote the projection on the section along 
the leaves of this foliation. Let $\Phi^{\sst J}_t$ be the Hamiltonian flow of $J$, 
with flow parameter $t$. 
\newline It may be remarked that, if a global transversal section is not at hand,
neighbourhoods on that individual operators $S$ can be constructed in the manner 
described below can be put together to yield a common $S$ operator on the union of the 
neighbourhoods, as long as their intersection fulfills certain requirements; since it is 
the aim of the present discussion to outline the main ideas of the star product reduction 
presented here, we do not embark on giving the details of this problem.

\begin{lemma} \label{S} \EinrL
Let $p \equiv \Phi^{\sst J}_{t(p)} ({\mathfrak p}(p))$ be any point on $\Cc$ and $f \in 
\CPl$, and let $S = \sum_{n=0}^{\infty} \la^n S_n$ be inductively defined as 
     \begin{eqnarray*}
          S_0 &=& \id \\ (S_nf)(p) &:=& -i \int_0^{t(p)} {\Phi^{\sst J}_t}^* (F_{n+1}[S_0, 
          \ldots , S_{n-1};T_0, \ldots ,T_n](f))({\mathfrak p}(p))\: dt \\ &&\mbox{\it for 
          } n \geq 1 \mbox{\it , where (for $n \geq 2$)} \\ && \textstyle F_n[S_0, \ldots 
          ,S_{n-2};T_0, \ldots ,T_{n-1}](f) := 
             \sum_{k=2}^{n} \M{k}(J,S_{n-k}f)  \\
          && \textstyle \hspace*{22.6ex} + \sum_{i=1}^{n-1} \sum_{k=1}^{n-i}
             T_i(\M{k}(J,S_{n-k-i}f))\:,
     \end{eqnarray*}
and in the sandwich neighborhood $U \supset \Cc$ we set $S_nf := \prol(S_nf)$.
\newline Then $S : B \cap F \rightarrow B^* \cap F$ is a linear isomorphism.
\end{lemma}             

{\bf Proof.} By construction we have for $n \geq 1$ and $f \in B \cap F$ that 
$L_{X_J}(S_nf)(p) = -iF_{n+1}(f)(p)$, $p \in \Cc$, so $-\{ J,S_nf \} = -iF_{n+1}(f)$ on 
$\Cc$. But a direct computation shows $\{J,S_nf\} -iF_{n+1}(f) = (T(J \ast Sf - Sf \ast 
J))_{n+1}$, that is $T(J  \ast Sf - Sf \ast J) = 0$ on $\Cc$ at order $n \geq 2$ with our 
$S$. At order $1$, $T(J \ast Sf - Sf \ast J) = i\{J,S_0f\} = 0$ on $\Cc$ as $f \in B$. All  
in all, $T(J \ast Sf - Sf \ast J) = h \cdot J$ for some $h \in \CPl$, because it is an 
element of $I$. Applying $T^{-1}$ to both sides and using iii) of lemma \ref{T} yields $J 
\ast Sf = g \ast J$ for some $g\in \CPl$, showing that $Sf \in B^*$ according to equation 
\ref{B^*}. $Sf \in F$ is clear by construction, and injectivity of $S : B \cap F 
\rightarrow B^* \cap F$ follows from $S_0 = \id$. So it remains to show that $S$ is onto.
To this end, we fix an arbitrary $f = \sum_{n=0}^{\infty} \la^n f_n \in B^* \cap F$. From 
this $f$, we construct a sequence $(g^{\sst (k)})_{k \in \Nl}$ of power series in $\la$, 
inductively defined by $g^{\sst (0)} := f$, ..., $g^{\sst (n+1)} := \frac{1}{\la}\, 
(g^{\sst (n)} - Sg^{\sst (n)}_0)$, and then in turn we can write down $g := g^{\sst (0)}_0 
+ \la g^{\sst (1)}_0 + \la^2 g^{\sst (2)}_0 + \cdots$, picking always the $0$th order term 
out of every series $g^{\sst (k)}$. It is not difficult to show now that $g$ is a 
well--defined power series in $\la$, that $g \in B \cap F$ and $Sg = f$. $\BBox$

For the locality of the future reduced star product, the following lemma is 
essential.

\begin{lemma} \label{suppS}
At every order $n$ and for every $f \in B \cap F$, $\supp{S_nf} \:\subseteq\: \supp{f}$.
\end{lemma}

{\bf Proof.} Because $f \in F$ and $S_nf \in F$, it is sufficient to consider the 
intersection of the supports with $\Cc$; so we have to prove that $\supp{S_nf} \cap \Cc 
\:\subseteq\: \supp{f} \cap \Cc$, $n \in \Nl$, $f \in B \cap F$. Suppose we already had proved 
this for $1,\ldots,n-1$. Then it follows from the construction of $S_n$ (using lemma 
\ref{T}, iv) in the step $\supp{F_n(f)} \cap \Cc \:\subseteq \supp{f} \cap \Cc$) with 
${\mathfrak p} : \Cc \rightarrow \sigma$, that ${\mathfrak p}(\supp{S_nf} \cap \Cc) 
\:\subseteq\: {\mathfrak p}(\supp{f} \cap \Cc)$. But $f \in B$, that is $f(p) = f({\mathfrak  
p}(p))$ for all $p 
\in \Cc$, so for every set $A \subseteq \Cc$, the implication ${\mathfrak p}(A) \:\subseteq\:
{\mathfrak p}(\supp{f} \cap \Cc) \Longrightarrow A \:\subseteq\: \supp{f} \cap \Cc$ holds. 
$\BBox$

\subsection{Construction of the Reduced Star Product}

\begin{lemma} \EinrL
The spaces $B^*/I^*$ and $B^* \cap F$ are linearly isomorphic through $\prol \circ T$. 
Likewise, $B/I$ and $B \cap F$ are isomorphic through $\prol$.
\end{lemma}

{\bf Proof.} Let $b \in B^*$, then $b = f+i$, $f \in F$, $i \in I^*$. But $i \in I^* 
\subset B^* \Longrightarrow f = b-i \in B^* \Longrightarrow f \in B^* \cap F$. So $b \in 
B^* \cap F \oplus I$. Conversely, $B^* \subset \CPl$ is a subspace, so trivially 
$B^* \cap F \oplus B^* \cap I^* \subset B^*$ and with $B^* \cap I^* = I^*$ holds $B^* \cap  
F \oplus I^* \subset B^*$. From $B^* = B^* \cap F \oplus I^*$ then, we see that (comparing 
this with the unique decomposition in equation \ref{decompT}) $\prol(T(f)) \in B^* \cap F$ 
for $f \in B^*$, and furthermore $\prol \circ T$ is well--defined on $B^*/I^*$, because 
for $f \in B^*$, $i \in I^*$, $\prol(T(f+i)) = \prol(T(f))$ since $T(i) \in I$. $\BBox$

We now have the following chain of linear isomorphisms:
     \be
          \CCdo \:\cong\: B/I \:\cong\: B \cap F \:\cong\: B^* \cap F \:\cong\: 
          (B^*/I^*,\ast) \:, 
     \ee 
the latter space endowed with a star product inherited from that of $(P, \omega)$. Now let 
$f,g \in B \cap F$, then $Sf, Sg \in B^* \cap F$. Regarding them as representants in $B^*/I^*$,
we form $Sg \ast Sf$ (this is in $B^*$, but not in $F$ anymore), so that $\prol(T(Sf \ast SG))
\in B^* \cap F$. Finally, applying $S^{-1}$ brings us back to $B \cap F$.

\begin{lemma} \label{RSP} \EinrL
and definition. Identifying $B \cap F$ and $\CCdo$, we set
     \begin{displaymath}
          f \star g := S^{-1}(\prol(T(Sf \ast Sg)))
     \end{displaymath}
for $f,g \in B \cap F$ and have a (local) star product on $B \cap F \:\cong\: \CCdo$. We
call it ``reduced from $\ast$'' by the first class constraint $J$.
\end{lemma}

{\bf Proof.} The properties of a star product (basically clear by construction) can be 
checked one by one, remembering $S_0 = (S^{-1})_0 = T_0 = \id$ and $\prol(f \cdot g) = 
\prol(f) \cdot \prol(g)$; to prove $f \star 1 = 1 \star f = f$, we use that $S$ vanishes on 
constants and $\prol(T(f)) = f$ for $f \in F$ (lemma \ref{T} ii)), so $f \star 1 = 
S^{-1}(\prol(T(Sf \ast 1))) = S^{-1}(\prol(T(Sf))) = S^{-1}(Sf) = f$ (since $Sf \in B^* 
\cap F \:\subset\: F$). Locality follows from lemma \ref{T} iv) and lemma \ref{suppS}. $\BBox$

If the $M_n$ are bidifferential operators of finite order, so are the $\tilde{M}{}_n$ 
associated with $\star$, as can be seen by confirming that neither $T$ nor $S$ can 
increase the number of derivatives.

\subsection{Example: Moyal product on $\Rl^{2n}$ to $\Rl^{2n-2}$}

We first show that in a most simple example, the reduction formalism gives the expected 
result. To this purpose, we take $\Rl^{2n}$ with the usual symplectic form and a global 
chart $(q^1,..,q^n;p_1,..,p_n)$. We feed the constraint function $J(\bom{q};\bom{p}) := 
p_n$ into the classical reduction formalism and get $\Rl^{2n-2}$ as Reduced Phase Space, 
for which $(q^1,..,q^{n-1};p_1,..,p_{n-1})$ may serve as a global chart. On $\Rl^{2n}$, we  
suppose the {\bf Moyal product} as given (for operator orderings see e.g. \cite{AW70}); 
its explicit form is $f \ast g = 
\sum_{r=0}^{\infty} \frac{1}{r!} (\frac{i \la}{2})^r\, \Lambda^{k_1 l_1} \ldots 
\Lambda^{k_r l_r} \frac{\partial^r f}{\partial \xi^{k_1} \ldots \partial \xi^{k_r}} 
\frac{\partial^r g}{\partial \xi^{l_1} \ldots \partial \xi^{l_r}}\:$, $l_s, k_s = 1 \ldots 
2n$, where $\Lambda$ denotes the Poisson tensor to $\omega = dq^i \wedge dp_i$.
\newline For the star product reduction, we have to choose a prolongation prescription 
off the constraint surface $\Cc = \{ (\bom{q};\bom{p}) \in \Rl^{2n} \mde p_n = 0 \}$, and 
we do this in the simplest manner by setting $(\prol (f))(\bom{q};p_1,..,p_n) := 
f(\bom{q};p_1,..,p_{n-1},0)$. In this case $\pj(f)$ is nothing else but a difference 
quotient in the direction of $p_n$. It turns out that the operator $T : I \rightarrow I^*$ 
can in this example be written as $T = \frac{1}{1-\la K}$ with $Kf := -M_1(\pj(f),J)$ 
being $\frac{1}{4}i$ times the difference quotient of $\frac{\partial f}{\partial q^n}$ in  
the direction of $p_n$, and a short calculation shows that $S : B \cap F \rightarrow B^* 
\cap F$ is equal to (the prolongation of) $\id - \la K$. But $B \cap F$ can be recognized as
the space of series of functions not depending on $q^n$ and $p_n$ (which is clear because 
$\Cinfty(\Rl^{2n-2})\pr \:\cong\: B/I \:\cong\: B \cap F$), so $S = \id : B \cap F 
\rightarrow B^* \cap F$; the equality of the spaces $B \cap F$ and $B^* \cap F$ can, of 
course, be established also in a direct way. But furthermore $T(f \ast g) = f \ast g$ for 
$f,g \in B \cap F$ on the basis of the Moyal product's special form, as well as $\prol(f 
\ast g) = f \ast g$, so by putting all this together we end with $f \star g = f \ast g$:
the Reduced Star Product is just again the Moyal product, this time for functions on 
$\Rl^{2n-2}$.

\subsection{Example: Wick product from $\Cx^{n+1}$ to $\CPn$}

Things look different if we reduce the Wick product (see for example again \cite{AW70})
from $\Cxwz$ to $\CPn$. 
This has already been done, even resulting in an explicit formula \cite{BBEW96a}, but taking 
into account the symmetries of the problem --- which we will ignore. We can therefore try 
our reduction mechanism on two different constraint functions, and we will obtain explicit  
formulae in both cases, for two star products on $\CPn$ that are not only different but 
inequivalent. 

In the first case, we consider $\Cxwz$ with usual symplectic form $\omega = 
\frac{i}{2}\, dz^i \wedge d\bar{z}{}^i$ and $J(z) = -\halb \zzq - \mu$, $\mu \in \Rl^{-}$ 
(where $\zzq$ abbreviates $\sum_{i=1}^{n+1} z^i \bar{z}{}^i$). This $J$ is an $\mbox{\it 
ad}_*$--equivariant momentum mapping for the $U(1)$ group action $z \mapsto e^{i\varphi}z$ 
on $\Cxwz$, but we agreed to ignore these aspects altogether.
We regard $J$ as first class constraint only 
and see that $\Cc := J^{-1}(0)$ is an immersed sphere $S^{2n+1} \imbed 
\Cxwz$ with radius $\sqrt{-2\mu}$ and that $\Cd$, the classical Reduced Phase Space, is 
just $\CPn$. On $\Cxwz$, let the {\bf Wick product} be given: $f \ast g = 
\sum_{r=0}^{\infty} \frac{2 \la^r}{r!} \sum_{i_1 = \cdots = i_r = 1}^{n} 
\frac{\partial^r f}{\partial z^{i_1} \ldots \partial z^{i_r}} \frac{\partial^r 
g}{\partial \bar{z}{}^{i_1} \ldots \partial \bar{z}^{i_r}}$.

We choose to prolongate every $f \in \Cinfty(\Cxwz)\pr$ off $S^{2n+1}$ in a radial way, 
that is we set $\prol(f)(z) := f(\p(z))$ where $\p :  \Cxwz \rightarrow S^{2n+1}$; $z 
\mapsto (\sqrt{-2\mu / \zzq}) z$ projects radially on $S^{2n+1}$. From the general 
formula \ref{pijott}, 
we get $\pj(f)(z) = \frac{f(\ssf{p}(z))- f(z)}{\frac{1}{2} \zzq + \mu}$. Onto $\Cc = S^{2n+1}$, 
this is continued as $\resC \pj(f) = \frac{1}{2\mu} (E + \Eq)f$, where $E$ and $\Eq$ 
denote the {\it Euler operators} $E = z^k \frac{\partial}{\partial z^k}$ and $\Eq = 
\bar{z}{}^k \frac{\partial}{\partial \bar{z}{}^k}$. Even in this example, the inductive 
formula \ref{Tn} for $T_n$ can be resolved in terms of the operator $K := \halb E \circ 
\pj$, yielding $T = \sum_{r=0}^{\infty} \la^r K^r = \frac{1}{1-\la K}$. It is important 
that $K$, like $\pj$, can be expressed by the Euler operators $E$ and $\Eq$ when 
evaluated on $\Cc$: it is $\resC (Kf) = \frac{1}{8\mu}(\halb(E^2-\Eq{}^2) + (E+\Eq) - 
(E+\Eq)E)f$. Another property of $K$ which will be of some importance later is that $K(f 
\cdot h) = K(f) \cdot h$ if $h \in \CPl$ is {\it homogeneous}, the latter meaning 
that $h(z) = h(\la z)$ for all $\la \in \Cx \!\setminus\! \{0\}$ or equivalently, that $h 
= \pi^* \eta$ with an $\eta \in \Cinfty(\CPn)\pr$ and $\pi : \Cxwz \rightarrow \CPn$ the 
canonical projection.

The space of homogeneous series is equal to $B \cap F$, as an analysis of the conditions 
$f \in F$ and $f \in B$ will show (resulting in $f \in B \cap F \Longleftrightarrow Ef = 
0$, $\Eq f = 0$) and like it should be because of $B \cap F \:\cong\: \CCdl$. The operator 
$S$ has to be evaluated on $B \cap F$. Its general recursive definition \ref{S} takes the 
form of $(E-\Eq)(S_nf) = -\sum_{k=1}^{n} K^k((E-\Eq)S_{n-k}f)$ in the present example, 
leading to $(E-\Eq)S_0 = E-\Eq$ and $(E-\Eq)S_1 = -K(E-\Eq)$ with higher orders vanishing. 
Because on $\Cc = S^{2n+1}$, $K$ can be expressed in terms of the Euler operators, $S = 
\id -\la K$ is a solution for $S$ on $\Cc$. But $K$ vanishes on $B \cap F$, so $S = \id : 
B \cap F \rightarrow B^* \cap F$. The equality of $B \cap F$ and its ``deformed'' 
counterpart $B^* \cap F$ can of course be established by direct computations also. After 
putting all this together, lemma \ref{RSP} gives us: let $f,g \in B \cap F$ be two 
homogeneous series of functions, $K$ as defined above and $\p(z) = (\sqrt{-2\mu / \zzq}) 
z$. Then $(f \star g)(z) := (\frac{1}{1-\la K}(f \ast g))(\p(z))$ is homogeneous and 
$\star$ thereby defines a star product on $\CPn$, reduced from the Wick product $\ast$ on 
$\Cxwz$. 

The Reduced Star Product, though, can be considerably simplified by the following 
considerations. We define now bidifferential operators $\tilde{M}{}_r$ from the Wick 
product operators $M_r$ by $\tilde{M}{}_r(f,g) := (\zzq)^r M_r(f,g)$ and observe that for 
$f,g$ homogeneous, $\tilde{M}{}_r(f,g)$ is again homogeneous. We already mentioned that 
$K$, applied to a product of which one factor is homogeneous, this factor can be passed 
through, so in $\frac{1}{1-\la K} \sum_{r=0}^{\infty}(\frac{\la}{\zzq})^r 
\tilde{M}{}_r(f,g)$, only terms of the form $K^l(\frac{1}{(\zzq)^k})(\p(z))$ remain to be 
evaluated. The result of the ensuing computations is the star product on $\CPn$, reduced 
from the Wick product $f \ast g = \sum_{k=0}^{\infty} \la^k M_k(f,g)$:
     \begin{eqnarray} \label{Sternprod1}
          f \star g &=& \sum_{k=0}^{\infty} \sum_{l=0}^{\infty} 
                        {\textstyle (\frac{\la}{-2\mu})^{k+l}}
                        A^{\s{(k)}}_l (\zzq)^k M_k(f,g) \nonumber \\
                    && A^{\s{(k)}}_l := (-1)^l \sum_{i_1 =1}^{k} \sum_{i_2 = 1}^{i_1}
                       \cdots \sum_{i_l =1}^{i_{l-1}}\, i_1 i_2 \cdots i_l
     \end{eqnarray} 
     
The numbers $A^{\s{(k)}}_l$ fulfil a variety of inductive relations that in turn can be used to 
gain another direct formula, namely $A^{\s{(k)}}_l = \frac{1}{(k-1)!} \sum_{n=1}^{k}
{k-1 \choose n-1} (-1)^{k+l-n} n^{k+l-1}$. It is $A^{\s{(1)}}_l = (-1)^l$, $A^{\s{(k)}}_0 = 1$
and $A^{\s{(k)}}_1 = -\halb k(k+1)$. This direct formula shows that these numbers are  
the same as the numbers just so called in \cite{BBEW96a} where they were obtained in an altogether
different way, thereby proving that the star product on $\CPn$ constructed in \cite{BBEW96a} and
ours are identical.

\subsection{Example: Wick product from $\Cx^{n+1}$ to $\CPn$, inequivalently}

But because our constraint function $J$ need not necessarily be a momentum mapping, we can                           
repeat the whole reduction process with a different $J$, for example $J(z) := \frac{1}{4} 
(\zzq)^2 - \mu^2$, $\mu \in \Rl^{-}$. The classical Reduced Phase Space is of course 
$\CPn$ in both cases, the constraint submanifold being the same. The Reduced Star 
Products, however, turn out to be not just different but inequivalent, as we will see. The 
calculations proceed along the lines already followed, so there is no need to go into the 
details. Let it be sufficient to mention that $T_n = \sum_{k=0}^{[\frac{n}{2}]} 
\frac{1}{k!} \frac{d^k}{d \nu^k} \left. (P + \nu R)^{n-k} 
\right |_{\!\nu=0}$ for $(Pf)(z) := -\halb 
\zzq E(\pj(f))(z)$ and $(Rf)(z) := -\frac{1}{4}\, E^2(\pj(f))(z)$ ($[\frac{n}{2}]$ 
denoting the integer part of $\frac{n}{2}$), that $B \cap F$ is -- of course! -- again the 
space of homogeneous series of functions (remark that we did not touch the prolongation 
description), that $S$ can be chosen as identity, and that both $P$ and $R$ do not ``see''
the homogeneous factors in their arguments. So, for $f$ and $g$ homogeneous, 
     \begin{eqnarray}
          f \star g &=& \sum_{k=0}^{\infty} \sum_{l=0}^{\infty} 
                        {\textstyle (\frac{\la}{-2\mu})^{k+l}}
                        B^{\s{(k)}}_l (\zzq)^k M_k(f,g) \nonumber \\
                    && B^{\s{(k)}}_l := (-2\mu)^{k+l} \sum_{j=0}^{[\frac{l}{2}]}
                       \frac{1}{j!} \frac{d^j}{d \nu^j}
                       ((-\halb \zzq E + {\textstyle \frac{1}{4}}\, \nu E^2) \circ 
                       \pj)^{l-j}\, \left. \frac{1}{(\zzq)^k} 
                       \right |_{\!{\nu=0 \atop \zzq = -2\mu}} \:,
     \end{eqnarray} 
is another star product on $\CPn$ reduced from the Wick product on $\Cxwz$.

Let us denote the product of formula \ref{Sternprod1} with $\tilde{\star}$, then we 
obtain by a straightforward computation $(f \star g - f \,\tilde{\star}\, g)_2 - (g \star f - 
g \,\tilde{\star}\, f)_2 = \halb i (\frac{\la}{-2\mu})^2 \zzq \{f,g\}$ to the second order. 
But $\halb i (\frac{\la}{-2\mu})^2 \zzq \{f,g\} = \halb \frac{\la}{-2\mu} (f \star g - g 
\star f)_1 = \frac{i}{2}\, \frac{\la^2}{-2\mu} \{f,g\}_{\CPn}$, where the Poisson bracket 
$\{ \cdot , \cdot \}_{\CPn}$ belongs to the symplectic Fubini-Study form 
$\omega_{\CPn}$ on $\CPn$, and $\omega_{\CPn}$ is not exact. On the other 
hand, a result in \cite{BCG97} (see also \cite{NT95a}) says that two equivalent star products, 
equal up to the 
order $k$, have necessarily an exact two--form as the antisymmetric part of their 
difference at order $k+1$. Because $(f \star g)_1 = (f \,\tilde{\star}\, g)_1$ in our 
examples, this theorem applies and we conclude that our two reduced star products cannot 
be equivalent. Roughly speaking, two different constraint functions, though inducing the 
same classical Reduced Phase Space, may lead to inequivalent quantum systems.

\subsection*{Acknowledgement}

The author is greatly indebted to M. Bordemann, who suggested to the author 
the ansatz to define the algebras $I^*$ and $B^*$ and accompanied the ensuing work critically.

\end{document}